\documentclass[11pt]{amsart}
\usepackage{amsmath}
\usepackage{amssymb}
\usepackage{amsfonts}
\usepackage{amsthm}
\usepackage{latexsym}
\usepackage{graphics}
\usepackage{graphicx}
\usepackage[all]{xy}
\usepackage{verbatim}

\newtheorem{Definition}{Definition}
\newtheorem{Lemma}{Lemma}

\newtheorem{Theorem}{Theorem}
\newtheorem{Proposition}{Proposition}
\newtheorem{Corollary}{Corollary}

\theoremstyle{definition}
\newtheorem*{Ack}{Acknowledgement}

\DeclareMathOperator{\End}{End}
\DeclareMathOperator{\sign}{sign}

\title{\bf Invariance of the BFV--complex}

\author[F. Sch\"atz]{Florian Sch\"atz}
\address{Institut f\"ur Mathematik, Universit\"at Z\"urich--Irchel, Winterthurerstrasse 190, CH-8057 Z\"urich, Switzerland}
\email{florian.schaetz@gmail.com}

\thanks{The author acknowledges partial support by the research grant of the University of Zurich,
by the SNF-grant 200020-121640/1,
by the European Union through the FP6 Marie Curie RTN ENIGMA (contract number MRTN-CT-2004-5652), 
and by the European Science Foundation through the MISGAM program.}

\begin{document}

\maketitle 

\begin{abstract}
The BFV-formalism was introduced to handle classical systems, equipped with symmetries. 
It associates a differential graded Poisson algebra
to any coisotropic submanifold $S$ of a Poisson manifold $(M,\Pi)$.

However the assignment (coisotropic submanifold) $\leadsto$ (differential graded Poisson algebra) is not canonical, since
in the construction several choices have to be made. One has to fix: 1. an embedding of the normal bundle $NS$ of $S$ into $M$ as a
tubular neighbourhood, 2. a connection $\nabla$ on $NS$ and 3. a special element
$\Omega$.

We show that different choices of a connection and an element $\Omega$
-- but with the tubular neighbourhood fixed -- lead to isomorphic differential graded Poisson algebras.
If the tubular neighbourhood is changed too, invariance
can be restored at the level of germs.
\end{abstract}

\section{Introduction}\label{s:intro}
The Batalin-Vilkovisky-Fradkin complex (BFV-complex for short) was introduced in order to understand physical systems with
complicated symmetries (\cite{BatalinFradkin}, \cite{BatalinVilkovisky}). The connection to homological algebra
was made explicit in \cite{Stasheff} later on. We focus on the smooth setting, i.e. we want to
consider arbitrary coisotropic
submanifolds of smooth finite dimensional Poisson manifolds. Bordemann and Herbig found a convenient adaptation
of the BFV-construction in this framework (\cite{Bordemann}, \cite{Herbig}): One obtains a differential graded
Poisson algebra associated to any coisotropic submanifold.
In \cite{Schaetz} a slight modification of the construction of Bordemann and Herbig was presented. It made
use of the language of higher homotopy structures and provided in particular a conceptual construction of the BFV-bracket.

Note that in the smooth setting the construction of the BFV-complex requires a choice of the following pieces of data:
1. an embedding of the normal bundle of the coisotropic submanifold as a tubular neighbourhood into the ambient Poisson manifold, 2. a connection
on the normal bundle, 3. a special function on a smooth graded manifold, called a BFV-charge.

We apply the point of view established in \cite{Schaetz} to clarify the dependence of the resulting
BFV-complex on these data.
If one leaves the embedding fixed and only changes the connection and the BFV-charge, one simply obtains two
isomorphic differential graded Poisson algebras, see Theorem \ref{theoremA} in Section \ref{s:connection}.
Note that the dependence on the choice of BFV-charge was well understood, see \cite{Stasheff} for instance.
Dependence on the embedding is more subtle. We introduce the notion of ``restriction'' of a given BFV-complex
to an open neighbourhood of the coisotropic submanifold inside its normal bundle (Definition \ref{d:restrictedBFV})
and show that different choices of embeddings lead to isomorphic restricted BFV-complexes -- see Theorem \ref{thm:central} in
Section \ref{s:tubularneighbourhood}. As a Corollary one obtains that a germ-version of the BFV-complex is independent of
all the choices up to isomorphism (Corollary \ref{cor:central}).

It turns out that the differential graded Poisson algebra associated to a fixed embedding of the normal bundle as a tubular neighbourhood,
yields a description of the moduli space of coisotropic sections in terms of the BFV-complex -- see \cite{Schaetz2}.

\begin{Ack}
I thank Alberto Cattaneo for remarks on a draft of this work. Moreover, I thank the referee for helpful comments.
\end{Ack}

\section{Preliminaries}\label{s:pre}
The purpose of this Section is threefold: to recollect some facts about the theory of higher homotopy structures,
to recall some concepts concerning Poisson manifolds and coisotropic submanifolds and to
outline the construction of the BFV-complex. More details on these subjects
can be found in Sections 2 and 3 of \cite{Schaetz} and in the references cited therein.
We assume the reader to be familiar with the theory of graded algebras and smooth graded manifolds.

\subsection{$L_{\infty}$-algebras: Homotopy Transfer and Homotopies}\label{ss:homotopies}
Let $V$ be a $\mathbb{Z}$-graded vector space over $\mathbb{R}$ (or any other field of characteristic $0$); i.e.,
$V$ is a collection $(V_{i})_{i \in \mathbb{Z}}$ of vector spaces $V_{i}$ over $\mathbb{R}$.
The homogeneous elements of $V$ of degree $i\in \mathbb{Z}$ are the elements of $V_{i}$. We denote the degree
of a homogeneous element $x \in V$ by $|x|$.  A morphism $f: V \to W$ of graded vector
spaces is a collection $(f_{i}: V_{i} \to W_{i})_{i \in \mathbb{Z}}$ of
linear maps.
The $n$th suspension functor $[n]$ from the category of graded vector spaces to itself is defined as follows:
given a graded vector space $V$,
$V[n]$ denotes the graded vector space corresponding to the collection $V[n]_{i}:=V_{n+i}$.
The $n$th suspension of a morphism $f: V \to W$ of graded vector spaces is given
by the collection $(f[n]_{i}:=f_{n+i}:V_{n+i} \to W_{n+i})_{i \in \mathbb{Z}}$.
The tensor product of two graded vector spaces $V$ and $W$ over $\mathbb{R}$ is the graded vector whose component
in degree $k$ is given by
\begin{align*}
(V\otimes W)_{k}:=\bigoplus_{r+s=k}V_{r}\otimes W_{s}.
\end{align*}
The denote this graded vector space by $V\otimes W$.

The structure of a {\em flat $L_{\infty}[1]$-algebra}
on $V$ is given by a family of multilinear maps $(\mu^{k}:V^{\otimes k} \to V[1])_{k\ge 1}$ that satisfies:
\begin{itemize}
\item[(1)] $\mu^k(\cdots \otimes a \otimes b \otimes \cdots)=(-1)^{|a||b|}\mu^k(\cdots \otimes b \otimes a \otimes \cdots)$
holds for all $k\ge 1$ and all homogeneous elements $a, b$ of $V$.
\item[(2)] The family of {\em Jacobiators} $(J^k)_{k\ge 1}$ defined by
\begin{multline*}
J^{k}(x_{1} \cdots x_{n}):=\\=
\sum_{r+s=k} \sum_{\sigma \in (r,s)-\text{shuffles}} \mspace{-36mu} \sign(\sigma)\, \mu^{s+1}(\mu^{r}(x_{\sigma(1)} \otimes \cdots \otimes x_{\sigma(r)})
\otimes x_{\sigma(r+1)} \otimes \cdots \otimes x_{\sigma(n)}) 
\end{multline*}
vanishes identically. Here $\sign(\cdot)$ is the Koszul sign, i.e. the representation of $\Sigma_n$
on $V^{\otimes n}$ induced by mapping the transposition $(2,1)$ to $a\otimes b \mapsto (-1)^{|a||b|}b\otimes a$.
Moreover $(r,s)$-shuffles are permutations $\sigma$ of $\{1,\dots,k=r+s\}$ such that $\sigma(1) < \cdots < \sigma(r)$ and
$\sigma(r+1) < \cdots < \sigma(k)$.
\end{itemize}
Since we are only going to consider flat $L_{\infty}[1]$-algebras we will suppress the adjective
``flat'' from now on. In this case the vanishing of the first Jacobiator implies that $\mu^1$ is a coboundary
operator.
We remark that an $L_{\infty}[1]$-algebra structure on $V$ is equivalent to the more traditional notion of an
$L_{\infty}$-algebra structure on $V[-1]$, see \cite{Operads} for instance.

Given an  $L_{\infty}$-algebra structure $(\mu^k)_{k\ge 1}$ on $V$, there is a distinguished subset of $V_1$ that contains
elements $v \in V_1$ satisfying the {\em Maurer-Cartan equation} (MC-equation for short)
\begin{align*}
\sum_{k\ge 1}\frac{1}{k!}\mu^k(v\otimes \cdots \otimes v)=0.
\end{align*}
This set is called the set of {\em Maurer-Cartan elements} (MC-elements for short) of $V$.

Let $V$ be equipped with an $L_{\infty}$-algebra structure such that the coboundary operator $\mu^1$ decomposes into $d + \delta$
with $d^2=0=\delta^2$ and $d\circ \delta + \delta \circ d =0$. i.e. $(V,d,\delta)$ is a double complex. Then -- under
mild convergence assumptions -- it is possible to construct an $L_{\infty}$-algebra structure on $H(V,d)$ that
is ``isomorphic up to homotopy'' to the original $L_{\infty}$-algebra structure on $V$ (\cite{GugenheimLambe}). More concretely, one has to
fix an embedding $i$ of $H(V,d)$ into $V$, a projection $pr$ from $V$ to $H(V,d)$ and a homotopy operator $h$ (of degree $-1$) which satisfies
\begin{align*}
d\circ h + h\circ d= id_V - i\circ pr.
\end{align*}
We will also impose the following side-conditions for the sake of simplicity: 1.)$h\circ h=0$, 2.)$pr\circ h=0$ and 3.)$h\circ i=0$.
Then explicit formulae for the structure maps for
an $L_{\infty}$-algebras on $H(V,d)$ can be written down. These are given in terms of rooted planar
trees, see \cite{Schaetz} for a review. We will explain the construction in more detail later on for the examples which are relevant for our purpose.

Furthermore one obtains $L_{\infty}$-morphisms between $H(V,d)$ and $V$ that induce inverse maps on cohomology. 
Such
$L_{\infty}$-morphisms are called {\em $L_{\infty}$ quasi-isomorphisms}.

Consider the differential graded algebra $(\Omega([0,1]),d_{DR},\wedge)$ of smooth forms on the interval $I:=[0,1]$.
The inclusions of a point $\{*\}$ as $0\le s \le 1$ induces a chain map $ev_s: (\Omega(I),d_{DR}) \to (\mathbb{R},0)$
that is a morphisms of algebras. 
Given any $L_{\infty}$-algebra structure on $V$ there is a natural $L_{\infty}$-algebra structure on
$V\otimes \Omega(I)$ defined by
\begin{align*}
\tilde{\mu}^1(v\otimes \alpha):=\mu^1(v)\otimes \alpha + (-1)^{|v|}v\otimes d_{DR}\alpha
\end{align*}
and
\begin{align*}
\tilde{\mu}^k((v_1\otimes \alpha_1)\otimes \cdots \otimes (v_{k}\otimes \alpha_k)):=
(-1)^{\#}\mu^{k}(v_1\otimes \cdots \otimes v_k)\otimes (\alpha_1\wedge \cdots \wedge \alpha_k)
\end{align*}
for $k\ge 2$. Here $\#$ denotes the sign one picks up by assigning $(-1)^{|v_{i+1}||\alpha_i|}$ to passing $\alpha_i$
from the left-hand side of $v_{i+1}$ to the right-hand side (and replacing $\alpha_{i+1}$ by $\alpha_i \wedge \alpha_{i+1}$).

Following \cite{Operads}, we call two morphisms $f$ and $g$ from an $L_{\infty}$-algebra $A$ to $B$ {\em homotopic}
if there exists an $L_{\infty}$-morphism $F$ from $A$ to $B\otimes \Omega(I)$ such that
\begin{itemize}
\item $(id\otimes ev_0)\circ F=f$ and
\item $(id\otimes ev_1)\circ F =g$ hold.
\end{itemize}
This defines an equivalence relation on
the set of $L_{\infty}$-morphisms from $A$ to $B$.

Let $F$ be an $L_{\infty}$-morphism from $A$ to $B\otimes \Omega(I)$.
Consequently $f_{s}:=ev_{s}\circ F$ is an $L_{\infty}$ morphism between $A$ and $B$ for any $s \in I$. Given
a MC-element $v$ in $A$ one obtains a one-parameter family of MC-elements
\begin{align*}
w_s:=\sum_{k\ge 1}\frac{1}{n!}(f_s)_{k}(v\otimes \cdots \otimes v)
\end{align*}
of $B$. Here $(f_s)_k$ denotes the $k$th Taylor component of $f_s$.

In the main body of this paper we are only interested in the following particular case:
$B$ is a differential graded
Lie algebra (i.e. only the first and second structure maps are non-vanishing). Denote the graded Lie bracket by
$[\cdot,\cdot]$. Furthermore we assume that the differential $D$ is given by the adjoint action of a degree $+1$
 element $\Gamma$ that
satisfies $[\Gamma,\Gamma]=0$. The MC-equation for an element $w$ of $(B,D=[\Gamma,\cdot],[\cdot,\cdot])$ reads
\begin{align*}
[\Gamma + w, \Gamma + w]=0.
\end{align*}
From the  one-parameter family of MC-elements $w_s$ in $B$
one obtains a one-parameter family of differential graded Lie algebras on $B$ by setting 
\begin{align*}
D_s(\cdot):=[\Gamma+w_s,\cdot]
\end{align*}
while leaving the bracket unchanged.

How are the differential graded Lie algebras $(B,D_{s},[\cdot,\cdot])$ related for different values of $s\in I$?
To answer this question we first apply the $L_{\infty}$ morphism $F: A \leadsto B\otimes \Omega(I)$ to $v$ and obtain a MC-element $w(t)+u(t)dt$ in $B\otimes \Omega(I)$.
It is straightforward to check that $w(s)=w_s$ for all $s \in I$. Moreover the MC-equation in $B\otimes \Omega(I)$
splits up into
\begin{align*}
[\Gamma + w(t),\Gamma + w(t)]=0
\end{align*}
and
\begin{align*}
\frac{d}{dt}w(t)= [u(t), \Gamma + w(t)].
\end{align*}
The second equation implies that whenever the adjoint action of $u(t)$ on $B$ can be integrated to a one-parameter
family of automorphisms $(U(t))_{t\in I}$, $U(s)$ establishes an automorphism of $(B,[\cdot,\cdot])$ that maps $\Gamma+w(0)$ to $\Gamma+w(s)$
(for any $s \in I$). Consequently:

\begin{Lemma}\label{l:isomorphism}
Let $A$  and $(B,[\Gamma,\cdot],[\cdot,\cdot])$ be differential graded Lie algebras,
$v$ a MC-element in $A$ and $F$ an $L_{\infty}$ morphism
from $A$ to $B\otimes \Omega(I)$ such that
\begin{align*}
\sum_{k\ge 1}\frac{1}{k!}F_k(v\otimes \cdots \otimes v)
\end{align*}
is well-defined in $B\otimes \Omega(I)$. Denote this element by $w(t)+u(t)dt$.
Furthermore the flow equation
\begin{align*}
X(0)=b, \quad \frac{d}{dt}|_{t=s}X(t)=[u(s),X(s)], s\in I
\end{align*}
is assumed to have a unique solution for arbitrary $b\in B$.

Then the one-parameter family $U(t)$
of automorphisms of $B$ that integrates the adjoint action by $u(t)$ maps $\Gamma + w(0)$ to
$\Gamma+w(t)$. In particular $U(s)$ is an isomorphims of differential graded Lie algebras
\begin{align*}
(B,[\Gamma+w(0),\cdot],[\cdot,\cdot]) \to (B,[\Gamma+w(s),\cdot],[\cdot,\cdot])
\end{align*}
for arbitrary $s\in I$.
\end{Lemma} 

\subsection{Coisotropic Submanifolds}\label{ss:coisotropic}

We essentially follow \cite{Weinstein}, where more details can be found.
Let $M$ be a smooth, finite dimensional manifold. 
The bivector field $\Pi$ on $M$ is {\em Poisson} if the binary operation $\{\cdot,\cdot\}$ on $\mathcal{C}^{\infty}(M)$ given by
$(f,g)\mapsto <\Pi,df\wedge dg>$ satisfies the {\em Jacobi identity}, i.e.
\begin{align*}
\{f,\{g,h\}\}=\{\{f,g\},h\}+\{g,\{f,h\}\}
\end{align*}
holds for all smooth functions $f$, $g$ and $h$. Here $<\negthickspace-,-\negthickspace>$ denotes the natural pairing
between $TM$ and $T^{*}M$. 
Alternatively one can consider the graded algebra $\mathcal{V}(M)$ of multivector fields 
on $M$ equipped with the Schouten-Nijenhuis bracket $[\cdot,\cdot]_{SN}$. A bivector field $\Pi$ is Poisson if and only if
$[\Pi,\Pi]_{SN}=0$.

Associated to any
Poisson bivector field $\Pi$ on $M$ there is a vector bundle morphism
$\Pi^{\#}: T^{*}M \to TM$ given by contraction.
Consider a submanifold $S$ of $M$. The annihilator $N^{*}S$ of $TS$ is a subbundle of
$T^{*}M$. This
subbundle fits into a short exact sequence of vector bundles:
$$
\xymatrix{
0 \ar[r] & N^{*}S \ar[r] & T^{*}M|_{S} \ar[r] & T^{*}S \ar[r] & 0}.
$$

\begin{Definition}\label{coisotropic}
A submanifold $S$ of a smooth, finite dimensional Poisson manifold $(M,\Pi)$
is called \textsf{coisotropic} if the restriction of $\Pi^{\#}$ to $N^{*}S$
has image in $TS$.
\end{Definition}

There is an equivalent characterization
of coisotropic submanifolds: define the vanishing ideal of $S$ by 
\begin{align*}
\mathcal{I}_{S}:=\{f\in \mathcal{C}^{\infty}(M): f|_{S}=0\}. 
\end{align*}
A submanifold $S$ is coisotropic if and only if $\mathcal{I}_{C}$ is a Lie subalgebra
of $(\mathcal{C}^{\infty}(M),\{\cdot,\cdot\})$.

\subsection{The BFV-Complex}\label{ss:BFV}

The BFV-complex was introduced by Batalin, Fradkin and Vilkovisky with application in physics in mind
(\cite{BatalinFradkin}, \cite{BatalinVilkovisky}). Later on Stasheff (\cite{Stasheff}) gave an interpretation of the BFV-complex in terms of homological algebra.
The construction we present below is explained with more details in \cite{Schaetz}. It uses a globalization
of the BFV-complex for arbitrary coisotropic submanifolds found by Bordemann and Herbig (\cite{Bordemann}, \cite{Herbig}).

Let $S$ be a coisotropic submanifold of a smooth, finite dimensional Poisson manifold $(M,\Pi)$.
We outline the construction a differential graded Poisson algebra, which we call a {\em BFV-complex for $S$ in $(M,\Pi)$}.
The construction depends on the choice of three pieces of data: 1. an embedding of the normal
bundle of $S$ into $M$ as a tubular neighbourhood, 2. a connection on $NS$ and 3. a special
smooth function, called the charge, on a smooth graded manifold.

Denote the normal bundle of $S$ inside $M$ by $E$. Consider the graded vector bundle
$E^*[1]\oplus E[-1] \to S$ over $S$ and let $\mathcal{E}^*[1]\oplus \mathcal{E}[-1]\to E$ be
the pull back of $E^*[1]\oplus E[-1] \to S$ along $E\to S$.

We define $BFV(E)$ to be the space of smooth functions on the graded manifold
which is represented by the graded vector bundle $\mathcal{E}^{*}[1]\oplus \mathcal{E}[-1]$ over $E$.
In terms of sections one has $BFV(E)=\Gamma(\bigwedge(\mathcal{E})\otimes\bigwedge(\mathcal{E^{*}}))$.
This algebra carries a bigrading given by 
\begin{align*}
BFV^{(p,q)}(E):=\Gamma(\wedge^{p}\mathcal{E}\otimes\wedge^{q}\mathcal{E^{*}}).
\end{align*}
In physical terminology $p$ / $q$ is referred to as the
\textsf{ghost degree} / \textsf{ghost-momentum degree} respectively. One defines
\begin{align*}
BFV^{k}(E):=\bigoplus_{p-q=k}BFV^{(p,q)}(E)
\end{align*}
and calls $k$ the \textsf{total degree} (in physical terminology this is the ``ghost number'').
%

The smooth graded manifold $\mathcal{E}^{*}[1]\oplus \mathcal{E}[-1]$ comes equipped with a Poisson
bivector field $G$ given by the natural fibre pairing between $\mathcal{E}$ and $\mathcal{E}^{*}$, i.e. it is defined to be the natural contraction on $\Gamma(\mathcal{E})\otimes \Gamma(\mathcal{E}^{*})$ and extended
to a graded skew-symmetric biderivation of $BFV(E)$.

\noindent
\hspace{0cm}{\bf Choice 1.} Embedding.
\newline
Fix an embedding $\psi: E \hookrightarrow M$ of the normal bundle of $S$ into $M$. Hence
the normal bundle $E$ inherits a Poisson bivector field which we also denote by $\Pi$. (Keep in mind
that $\Pi$ depends on $\psi$!)

\noindent
\hspace{0cm}{\bf Choice 2.} Connection.
\newline
Next choose a connection on the vector bundle $E\to S$. This induces a connection on $\wedge E \otimes \wedge E^{*} \to S$ and via
pull back one obtains a connection $\nabla$ on $\wedge \mathcal{E} \otimes \wedge \mathcal{E}^{*} \to E$. We denote the corresponding horizontal
lift of multivector fields by
\begin{align*}
\iota_{\nabla}: \mathcal{V}(E) \to \mathcal{V}(\mathcal{E}^{*}[1] \otimes \mathcal{E}[-1]).
\end{align*}
It extends to an isomorphism of graded commutative unital associative algebras
\begin{align*}
\varphi: \mathcal{A}:=\mathcal{C}^{\infty}(T^{*}[1]E\oplus \mathcal{E}^{*}[1]\oplus \mathcal {E}[-1] \oplus \mathcal{E}[0] \oplus \mathcal{E}^{*}[2]) \to \mathcal{V}(\mathcal{E}^{*}[1]\oplus \mathcal{E}[-1]).
\end{align*}
Using $\varphi$ we lift $\Pi$ to a bivector field on $\mathcal{E}^{*}[1]\oplus \mathcal{E}[-1]$. Since $\varphi$ fails 
in general to be a morphism
of Gerstenhaber algebras, $\varphi(\Pi)$ is not a Poisson bivector field. Similarly the sum $G+\varphi(\Pi)$ fails
to be a Poisson bivector field in general. However the following Proposition provides an appropriate correction term:

\begin{Proposition}\label{prop:lift}
Let $\mathcal{E}$ be a finite rank vector bundle with connection $\nabla$ over a smooth, finite dimensional manifold $E$.
Consider the smooth graded manifold $\mathcal{E}^ {*}[1]\oplus \mathcal{E}[-1] \to E$ and denote the Poisson
bivector field on it coming from the natural fibre pairing between $\mathcal{E}$ and $\mathcal{E}^{*}$ by $G$.

Then there is an $L_{\infty}$ quasi-isomorphism $\mathcal{L}_{\nabla}$ between the graded Lie algebra
\begin{align*} 
(\mathcal{V}(E)[1],[\cdot,\cdot]_{SN})
\end{align*}
 and the differential graded Lie algebra 
\begin{align*}
(\mathcal{V}(\mathcal{E}^{*}[1]\oplus \mathcal{E}[-1])[1],[G,\cdot]_{SN},[\cdot,\cdot]_{SN}).
\end{align*}
\end{Proposition}

A proof of Proposition \ref{prop:lift} can be found in \cite{Schaetz}.
It immediately implies

\begin{Corollary}\label{cor:lift}
Let $\mathcal{E } \to E$ be a finite rank vector bundle with connection $\nabla$ over a smooth, finite dimensional Poisson manifold
$(E,\Pi)$. Consider the smooth graded manifold $\mathcal{E}^{*}[1]\oplus \mathcal{E}[-1] \to E$ and denote the Poisson
bivector field on it coming from the natural fibre pairing between $\mathcal{E}$ and $\mathcal{E}^{*}$ by $G$.

Then there is a Poisson bivector field $\hat{\Pi}$ on $\mathcal{E}^{*}[1]\oplus \mathcal{E}[-1]$ such that
\newline $\hat{\Pi}= G + \varphi(\Pi) + \triangle$ for $\triangle \in \mathcal{V}^{(1,1)}(\mathcal{E}^{*}[1]\oplus \mathcal{E}[-1])$. 
\end{Corollary}

For a proof we refer the reader to \cite{Schaetz} again.

We remark that $\mathcal{V}^{(1,1)}(\mathcal{E}^{*}[1]\oplus \mathcal{E}[-1])$ is the ideal of 
$\mathcal{V}(\mathcal{E}^{*}[1]\oplus \mathcal{E}[-1])$
generated by multiderivations which map any tensor product of functions of total bidegree 
$(p,q)$ to a function of bidegree $(P,Q)$ where $P>p$ and $Q>q$.
In general, let $\mathcal{V}^{(r,s)}(\mathcal{E}^{*}[1]\oplus \mathcal{E}[-1])$ be the ideal
generated by multiderivations of $\mathcal{C}^{\infty}(\mathcal{E}^{*}[1]\oplus \mathcal{E}[-1])$
with total ghost degree larger than or equal to $r$ and total ghost-momentum degree larger than or equal to $s$, respectively. 

The bivector field $\hat{\Pi}$ from Corollary \ref{cor:lift} equips $\mathcal{E}^{*}[1]\oplus \mathcal{E}[-1]$
with the structure of a graded Poisson manifold. Consequently $BFV(E)$ inherits a graded Poisson bracket which we
denote by $[\cdot,\cdot]_{BFV}$. It is called the {\em BFV-bracket}. Keep in mind that the BFV-bracket depends on the connection
on $E \to S$ we have chosen.

\noindent
\hspace{0cm}{\bf Choice 3.} Charge.
\newline
The last step in the construction of the BFV-complex is to provide a special solution to the MC-equation associated to
$(BFV(E),[\cdot,\cdot]_{BFV})$, i.e. one constructs a degree $+1$ element $\Omega$ that satisfies
\begin{align*}
[\Omega,\Omega]_{BFV}=0.
\end{align*}
Additionally, one requires this element $\Omega$ to contain the tautological section of $\mathcal{E}\to E$ as the lowest order term.
To be more precise, recall that 
\begin{align*}
BFV^{1}(E)=\bigoplus_{k\ge 0}\Gamma(\wedge^{k}\mathcal{E}\otimes \wedge^{k-1}\mathcal{E}^{*}).
\end{align*}
Hence any element of $BFV^{1}(E)$ contains a (possibly zero) component in $\Gamma(\mathcal{E})$. One requires
that the component of $\Omega$ in $\Gamma(\mathcal{E})$ is given by the tautological section of
$\mathcal{E}\to E$. A MC-element satisfying this requirement is called a {\em BFV-charge}.

\begin{Proposition}\label{prop:charge}
Let $(E,\Pi)$ be a vector bundle equipped with a Poisson bivector field and denote its zero section by $S$.
Fix a connection on $E\to S$ and equip the ghost/ghost-momentum bundle $\mathcal{E}^*[1]\oplus \mathcal{E}[-1] \to E$ with
the corresponding BFV-bracket $[\cdot,\cdot]_{BFV}$.

\begin{enumerate}
\item There is a degree $+1$ element $\Omega$ of $BFV(E)$ whose component in $\Gamma(\mathcal{E})$ is given by the tautological
section $\Omega_0$ and that satisfies
\begin{align*}
[\Omega,\Omega]_{BFV}=0
\end{align*}
if and only if $S$ is a coisotropic submanifold of $(E,\Pi)$.
\item Let $\Omega$ and $\Omega'$ be two BFV-charges. Then there is an automorphism
of the graded Poisson algebra $(BFV(E),[\cdot,\cdot]_{BFV})$ that maps $\Omega$ to $\Omega'$. 
\end{enumerate}
\end{Proposition}

See \cite{Stasheff} for a proof of this proposition.

Given a BFV-charge $\Omega$ one can define a differential $D_{BFV}(\cdot):=[\Omega,\cdot]_{BFV}$, called {\em BFV-differential}.
It is well-known that the cohomology with respect to $D$ is isomorphic to the Lie algebroid cohomology of $S$ (as a coisotropic
submanifold of $(E,\Pi)$).

By the second part of Proposition \ref{prop:charge}, different choices of the BFV-charge lead to isomorphic differential graded Poisson algebra
structures on $BFV(E)$. In the next Section we will establish that different choices of connection on $E\to S$ lead to differential Poisson
algebras that lie in the same isomorphism class. The dependence on the embedding of the normal bundle of $S$ is more subtle and will
be clarified in Section \ref{s:tubularneighbourhood}. 

\section{Choice of Connection}\label{s:connection}
Consider a vector bundle $E$ equipped with a Poisson bivector field $\Pi$ such that
that zero section $S$ is coisotropic.
The aim of this Section is to investigate the dependence of the differential graded 
Poisson algebra $(BFV(E),D_{BFV},[\cdot,\cdot]_{BFV})$ constructed in Subection \ref{ss:BFV} on the choice of a connection 
$\nabla$ on $E \to S$.

Recall that in order to lift the Poisson bivector field $\Pi$ to a bivector field on 
$\mathcal{E}^{*}[1]\oplus \mathcal{E}[-1]$, a connection $\nabla$ on $E \to S$ was used. 
Furthermore the 
$L_{\infty}$ quasi-isomorphism
between $(\mathcal{V}(E)[1],[\cdot,\cdot]_{SN})$ and 
$(\mathcal{V}(\mathcal{E}^{*}[1]\oplus \mathcal{E}[-1])[1],[G,\cdot]_{SN},[\cdot,\cdot]_{SN})$ in Proposition \ref{prop:lift} 
depends on $\nabla$ too.
Consequently so does the graded Poisson bracket $[\cdot,\cdot]_{BFV}$.

Let $\nabla_{0}$ and $\nabla_{1}$ be two connections on a smooth finite rank vector bundle
$\mathcal{E} \to E$. By Proposition \ref{prop:lift} we obtain two $L_{\infty}$ quasi-isomorphisms $\mathcal{L}_{\nabla_{0}}$ and 
$\mathcal{L}_{\nabla_{1}}$ from $(\mathcal{V}(E)[1],[\cdot,\cdot]_{SN})$ to
$(\mathcal{V}(\mathcal{E}^{*}[1]\oplus \mathcal{E}[-1])[1],[G,\cdot]_{SN},[\cdot,\cdot]_{SN})$.
Although these morphisms depend on the connections, this dependence is
very well-controlled:

\begin{Proposition}\label{prop:lift_homotopy}
Let $\mathcal{E}$ be a smooth finite rank vector bundle over a smooth, finite dimensional manifold $E$
equipped with two connections $\nabla_{0}$ and $\nabla_{1}$. Denote
the associated $L_{\infty}$ quasi-isomorphisms between  $(\mathcal{V}(E)[1],[\cdot,\cdot]_{SN})$
and $(\mathcal{V}(\mathcal{E}^{*}[1]\oplus \mathcal{E}[-1]),[G,\cdot]_{SN},[\cdot,\cdot]_{SN})$
from Proposition \ref{prop:lift} by $\mathcal{L}_{0}$ and $\mathcal{L}_{1}$ respectively.

Then there is an $L_{\infty}$ quasi-isomorphism
\begin{align*}
\hat{\mathcal{L}}: (\mathcal{V}(E)[1],[\cdot,\cdot]_{SN})\leadsto (\mathcal{V}(\mathcal{E}^{*}[1]\oplus \mathcal{E}[-1])\otimes \Omega(I),[G,\cdot]_{SN}+d_{DR},[\cdot,\cdot]_{SN})
\end{align*}
such that $(id\otimes ev_{0})\circ \hat{\mathcal{L}} = \mathcal{L}_{0}$ and
$(id\otimes ev_{1})\circ \hat{\mathcal{L}}=\mathcal{L}_{1}$ hold.
\end{Proposition}

\begin{proof}
Given two connections $\nabla_{0}$ and $\nabla_{1}$, one can define
a family of connections $\nabla_{s}:=\nabla_{0} + s(\nabla_{1} - \nabla_{0})$ parametrized by
the closed unit interval $I$. Consequently we obtain a one-parameter family of isomorphisms
of graded algebras
\begin{align*}
\varphi_{s}: \mathcal{A}:=\mathcal{C}^{\infty}(T^{*}[1]E\oplus \mathcal{E}^{*}[1]\oplus \mathcal{E}[-1] \oplus \mathcal{E}[0] \oplus \mathcal{E}^{*}[2]) \xrightarrow{\cong} \mathcal{V}(\mathcal{E}^{*}[1]\oplus \mathcal{E}[-1]),
\end{align*}
extending the horizontal lifting with respect to the connection $\nabla_{s}\oplus\nabla_{s}^{*}$.
Via this identification, $\mathcal{A}$
inherits a one-parameter family of Gerstenhaber brackets which we denote by $[\cdot,\cdot]_{s}$.
and a differential $\tilde{Q}$ which can be checked to be independet from $s$ in local coordinates.

For arbitrary $s \in I$ these structures fit into the following commutative
diagram:
\begin{align*}
\xymatrix{
 & (\mathcal{A}[1],\tilde{Q},[\cdot,\cdot]_{0}) \ar[dl]_{Pr} \ar[dr]^{\varphi_{0}} & &\\
(\mathcal{V}(E)[1],[\cdot,\cdot]_{SN}) & & (\mathcal{V}(\mathcal{E}^{*}[1]\oplus \mathcal{E}[-1])[1],[G,\cdot]_{SN},[\cdot,\cdot]_{SN})\\
& (\mathcal{A}[1],\tilde{Q},[\cdot,\cdot]_{s}) \ar[ul]^{Pr} \ar[ur]_{\varphi_{s}} \ar[uu]^{\psi_{s}} & &}
\end{align*}
where $\psi_{s}:=\varphi_{0}^{-1}\circ \varphi_{s}$ is a morphism of differential graded algebras and of
Gerstenhaber algebras. $Pr$ denotes the natural projection.

It is straightforward to show that the cohomology of $(\mathcal{A},\tilde{Q})$ is $\mathcal{V}(E)$ and that
the induced $L_{\infty}$ algebra coincides with $(\mathcal{V}(E)[1],[\cdot,\cdot])$, see the proof of Proposition 1 in \cite{Schaetz}.
Hence we obtain a one-parameter family
of $L_{\infty}$ quasi-isomorphisms $\mathcal{J}_s: (\mathcal{V}(E)[1],[\cdot,\cdot]_{SN}) \leadsto (\mathcal{A}[1],\tilde{Q},[\cdot,\cdot]_{s})$. 
Composition with $\psi_{s}$ yields a one-parameter family of $L_{\infty}$ quasi-isomorphisms
\begin{align*}
\mathcal{K}_{s}: (\mathcal{V}(E)[1],[\cdot,\cdot]_{SN}) \leadsto (\mathcal{A}[1],\tilde{Q},[\cdot,\cdot]_{0}).
\end{align*}
We remark that the composition of $\mathcal{J}_s$ with $\varphi_s$ yields the $L_{\infty}$ quasi-isomorphism $\mathcal{L}_s$ between
$(\mathcal{V}(E),[\cdot,\cdot]_{SN})$ and $(\mathcal{V}(\mathcal{E}^{*}[1]\oplus E[-1]),[G,\cdot]_{SN},[\cdot,\cdot]_{SN})$
associated to the connection $\nabla_{s}$ from Proposition \ref{prop:lift}.
Consequently $\mathcal{L}_{0}$ ($\mathcal{L}_{1}$) is the composition of $\mathcal{K}_{0}$ ($\mathcal{K}_{1}$)
with $\varphi_{0}$.

Next, consider the differential graded Lie algebra $(\mathcal{A}[1]\otimes \Omega(I), \tilde{Q}+d_{DR},[\cdot,\cdot]_{0})$.
To prove Proposition \ref{prop:lift}, a homotopy $\tilde{H}$ for $\tilde{Q}$ was constructed in \cite{Schaetz} such that
\begin{align*}
\tilde{Q}\circ \tilde{H} + \tilde{H}\circ \tilde{Q}=id - \iota \circ Pr
\end{align*}
is satisfied. Here, $\iota$ denotes the natural inclusion $\mathcal{V}(E)\hookrightarrow \mathcal{A}$.
One defines a one-parameter family of homotopies $\tilde{H}_{s}:=\psi_{s} \circ \tilde{H}\circ \psi_{s}^{-1}$ 
and checks that 
\begin{align*}
\tilde{Q}\circ \tilde{H}_{s}+\tilde{H}_{s}\circ \tilde{Q}=id - \psi_{s} \circ \iota \circ Pr 
\end{align*}
holds.

We define $\hat{Pr}: \mathcal{A}\otimes \Omega(I) \to \mathcal{V}(E)\otimes \Omega(I)$ to be $Pr\otimes id$
and $\hat{\iota}: \mathcal{V}(E)\otimes \Omega(I) \to \mathcal{A}\otimes\Omega(I)$
to be $\hat{\iota}:=(\psi_{s}\circ \iota)\otimes id$. Clearly $\hat{Pr}\circ \hat{\iota}=id$ and $\tilde{H}_{s}$
provides a homotopy between
$id$ and $\hat{\iota} \circ \hat{Pr}$.
Moreover the side-conditions $\tilde{H}_{s}\circ \tilde{H}_{s}=0$, $\hat{Pr}\circ \tilde{H}_{s}=0$
and $\tilde{H}_{s} \circ \hat{\iota} =0$ are still satisfied. We summarize the situation in the following
diagram:
\begin{align*}
\xymatrix{(\mathcal{V}(E)\otimes \Omega(I),0) \ar@<2pt>[r]^<<<<<{\hat{\iota}_s}& \ar@<2pt>[l]^>>>>>{\hat{Pr}} (\mathcal{A}\otimes \Omega(I),\tilde{Q})},\tilde{H}_{s}.
\end{align*}

Following Subsection \ref{ss:homotopies} these data can be used to perform homological transfer.
The input consists of the differential graded Lie algebra
\begin{align*}
(\mathcal{A}[1]\otimes \Omega(I),\tilde{Q}+d_{DR},[\cdot,\cdot]_{0}).
\end{align*}
To construct the induced structure maps,
one has to consider oriented rooted trees with bivalent and trivalent interior vertices. 
The leaves (the exterior vertices with the root excluded) are decorated
by $\hat{\iota}$, the root by $\hat{Pr}$, the interior bivalent vertices by $d_{DR}$, the interior
trivalent vertices by $[\cdot,\cdot]_{0}$ and the interior edges (i.e. the edges not connected to any exterior vertices)
by $-\tilde{H}_{s}$. One then composes these maps in the order given by the orientation towards the root.
The associated $L_{\infty}$ quasi-isomorphism is constructed in the same manner, however, the root
is not decorated by $\hat{Pr}$ but by $-\tilde{H}_{s}$ instead.

Recall that $\mathcal{V}^{(r,s)}(\mathcal{E}^{*}[1]\oplus \mathcal{E}[-1])$ is the ideal
generated by multiderivations of $\mathcal{C}^{\infty}(\mathcal{E}^{*}[1]\oplus \mathcal{E}[-1])$
with total ghost degree larger than or equal to $r$ and total ghost-momentum degree larger than or equal to $s$, respectively.
One can check inductively that trees decorated with $e$ copies of $-\tilde{H}_s$
increase the filtration index by $(e,e)$. Moreover
trees containing more than one interior bivalent vertex do not contribute since $d_{DR}$
increases the form-degree by $1$. These facts imply that 1. the induced structure is given by
$(\mathcal{V}(E)[1]\otimes \Omega(I),d_{DR},[\cdot,\cdot]_{SN})$ and 2. there is an $L_{\infty}$ quasi-isomorphism
\begin{align*}
(\mathcal{V}(E)[1]\otimes \Omega(I),d_{DR},[\cdot,\cdot]_{SN})\leadsto (\mathcal{A}[1]\otimes\Omega(I),\tilde{Q}+d_{DR},[-,-]_{0}). 
\end{align*}
We define 
\begin{align*}
\tilde{\mathcal{K}}: (\mathcal{V}(E)[1],[\cdot,\cdot]_{SN}) \leadsto (\mathcal{V}(\mathcal{A}[1]\otimes \Omega(I),\tilde{Q}+d_{DR},[\cdot,\cdot]_{SN})
\end{align*}
to be the composition
of this $L_{\infty}$ quasi-isomorphism and the obvious $L_{\infty}$ quasi-isomorphism
$(\mathcal{V}(E)[1],[\cdot,\cdot]_{SN})\hookrightarrow (\mathcal{V}(E)[1]\otimes \Omega(I),d_{DR},[\cdot,\cdot]_{SN})$.

The composition of $\hat{\mathcal{K}}$ with $id\otimes ev_{s}: \mathcal{A}\otimes \Omega(I) \to \mathcal{A}$
can be computed as follows: first of all only trees without any bivalent interior edges contribute since
all elements of form-degree $1$ vanish under $id \otimes ev_{s}$. Using the identities
$\psi_{s}^{-1}([\psi_{s}(-),\psi_{s}(-)]_{0})=[-,-]_{s}$,
$\tilde{H}_{s}=\psi_{s}\circ \tilde{H}\circ \psi_{s}^{-1}$ and $\hat{\iota}=\psi_{s}\circ \iota$ it is a straightforward
to show that $(id\otimes ev_{s})\circ \hat{\mathcal{K}}=\psi_{s}\circ\mathcal{K}_{s}$. Hence
\begin{align*}
\varphi_{0}\circ (id\otimes ev_{s})\circ \hat{\mathcal{K}}=\varphi_{s}\circ\mathcal{K}_{s}=\mathcal{L}_{s}.
\end{align*}
Finally, we define the $L_{\infty}$ quasi-isomorphism $\hat{\mathcal{L}}$ between $(\mathcal{V}(E)[1],[\cdot,\cdot]_{SN})$
and $(\mathcal{V}(\mathcal{E}^{*}[1]\oplus \mathcal{E}[-1])[1]\otimes \Omega(I),[G,\cdot]_{SN}+d_{DR},[\cdot,\cdot]_{SN})$
to be $(\varphi_{0}\otimes id) \circ \hat{\mathcal{K}}$. By construction $(id\otimes ev_0) \circ \hat{\mathcal{L}}=\mathcal{L}_0$ and
$(id\otimes ev_1)\circ \hat{\mathcal{L}}=\mathcal{L}_1$ are satisfied.
\end{proof}

We remark that Propositions \ref{prop:lift} and \ref{prop:lift_homotopy} seem to permit ``higher analogous'',
where one incorporates the differential graded algebra of differential forms on the n-simplex
$\Omega(\triangle^{n})$ instead of just $\Omega(\{*\})=\mathbb{R}$ (Proposition \ref{prop:lift})
or $\Omega(I)$ (Proposition \ref{prop:lift_homotopy}) -- see \cite{Costello}, where this idea was worked out in the
context of the BV-formalism.

\begin{Corollary}\label{cor:lift_homotopy}
Let $\mathcal{E}$ be a finite rank vector bundle over a smooth, finite dimensional Poisson manifold $(E,\Pi)$.
Suppose $\nabla_{0}$ and $\nabla_{1}$ are two connections on $\mathcal{E} \to E$. Denote
the associated $L_{\infty}$ quasi-isomorphisms between  $(\mathcal{V}(E)[1],[\cdot,\cdot]_{SN})$
and $(\mathcal{V}(\mathcal{E}^{*}[1]\oplus \mathcal{E}[-1])[1],[G,\cdot]_{SN},[\cdot,\cdot]_{SN})$
from Proposition \ref{prop:lift} by $\mathcal{L}_{0}$ and $\mathcal{L}_{1}$, respectively.
Applying these $L_{\infty}$ quasi-isomorphisms to $\Pi$ yields
two MC-elements $\tilde{\Pi}_{0}$ and $\tilde{\Pi}_{1}$ of 
$(\mathcal{V}(\mathcal{E}^{*}[1]\oplus \mathcal{E}[-1])[1],[G,\cdot]_{SN},[\cdot,\cdot]_{SN})$. 
Hence $\hat{\Pi}_0:=G+\tilde{\Pi}_{0}$ and $\hat{\Pi}_1:=G+\tilde{\Pi}_1$ are MC-elements of 
$(\mathcal{V}(\mathcal{E}^{*}[1]\oplus \mathcal{E}[-1])[1],[\cdot,\cdot]_{SN})$, i.e. Poisson bivector fields
on $\mathcal{E}^{*}[1]\oplus \mathcal{E}[-1]$.

There is a diffeomorphism of the smooth graded manifold $\mathcal{E}^{*}[1]\oplus \mathcal{E}[1]$ such that the
induced automorphism of $\mathcal{V}(\mathcal{E}^{*}[1]\oplus \mathcal{E}[-1])$ maps $\hat{\Pi}_{0}$ to $\hat{\Pi}_{1}$.
Moreover, this diffeomorphism induces a diffeomorphism of the base $E$ which coincides with the identity.
\end{Corollary}

\begin{proof}
Apply the $L_{\infty}$ quasi-isomorphism $\hat{\mathcal{L}}$ from Proposition \ref{prop:lift_homotopy} to $\Pi$ and add $G$ to obtain a MC-element $\hat{\Pi}+\hat{Z}dt$
of $(\mathcal{V}(\mathcal{E}^{*}[1]\oplus \mathcal{E}[-1])[1]\otimes \Omega(I),d_{DR},[\cdot,\cdot]_{SN})$. 
Let $\mathcal{L}_{s}$ denote the $L_{\infty}$ quasi-isomorphism from Proposition \ref{prop:lift}
constructed with the help of the connection $\nabla_{0}+s(\nabla_{1}-\nabla_{0})$. Recall that $(id\otimes ev_s)\circ \hat{\mathcal{L}}=\mathcal{L}_s$
holds for all $s\in I$.

We set $\hat{\Pi}_{s}:=(id\otimes ev_{s})(\hat{\Pi})$
and $\hat{Z}_{s}:=(id\otimes ev_{s})(\hat{Z})$. Proposition \ref{prop:lift_homotopy} implies that this definition of
$\hat{\Pi}_{s}$ is compatible with $\hat{\Pi}_{0}$ and $\hat{\Pi}_{1}$ defined in the Corollary.

We want to apply Lemma \ref{l:isomorphism} to $A:=(\mathcal{V}(E)[1],[\cdot,\cdot]_{SN})$,
$B:=(\mathcal{V}(\mathcal{E}^{*}[1]\oplus \mathcal{E}[-1])[1],[G,\cdot]_{SN},[\cdot,\cdot]_{SN})$ and $F:=\hat{\mathcal{L}}$.
To do so, it remains to show that the flow of $\hat{Z}_{s}$ is globally well-defined for $s \in [0,1]$. Recall that $\hat{Z}$
is the one-form part of the MC-element constructed from the Poisson
bivector field $\Pi$ on $E$ with help of the $L_{\infty}$ quasi-isomorphism
$\hat{\mathcal{L}}: (\mathcal{V}(E)[1],[\cdot,\cdot]_{SN}) \leadsto 
(\mathcal{V}(\mathcal{E}^{*}[1]\oplus \mathcal{E}[-1])\otimes \Omega(I),[G,\cdot]_{SN}+d_{DR},[\cdot,\cdot]_{SN})$. 
Only trees with exactly one bivalent interior vertex give non-zero
contributions because the form degree must be one. Consequently there is at least one homotopy in
the diagram and by the degree estimate in the proof of Proposition \ref{prop:lift_homotopy} this implies
that $\hat{Z}$ is contained in $\mathcal{V}^{(1,1)}(\mathcal{E}^{*}[1]\oplus \mathcal{E}[-1])\otimes\Omega(I)$. Hence
the derivation $[\hat{Z},-]_{SN}$ is nilpotent and can be integrated. Furthermore the degree estimate directly implies the last claim of the Corollary.
\end{proof}

The following is an immediate consequence of the previous Corollary:

\begin{Corollary}\label{cor:BFV-bracket_homotopy}
Let $(E,\Pi)$ be a vector bundle $E\to S$ equipped with a Poisson structure $\Pi$ such that $S$ is a coisotropic submanifold.
Fix two connections $\nabla_{0}$ and $\nabla_{1}$ on $E \to S$ and denote the corresponding
graded Poisson brackets on $BFV(E)$ by $[\cdot,\cdot]^{0}_{BFV}$
and $[\cdot,\cdot]^{1}_{BFV}$ respectively.

There is an isomorphism of graded Poisson algebra 
\begin{align*}
(BFV(E),[\cdot,\cdot]^{0}_{BFV}) \xrightarrow{\cong} (BFV(E),[\cdot,\cdot]^{0}_{BFV}).
\end{align*}
Moreover
the induced automorphism of $\mathcal{C}^{\infty}(E)$ coincides with the identity.
\end{Corollary}

Combining Proposition \ref{prop:charge} and Corollary \ref{cor:BFV-bracket_homotopy} we obtain

\begin{Theorem}\label{theoremA}
Let $E$ be a vector bundle equipped with a Poisson bivector $\Pi$ such that the zero Section $S$ is a coisotropic submanifold.
Recall that the pull back of $E\to S$ by $E\to S$ is denoted by $\mathcal{E}\to E$ 
and 
\begin{align*}
BFV(E):=\mathcal{C}^{\infty}(\mathcal{E}^{*}[1]\oplus \mathcal{E}[-1])=\Gamma(\wedge\mathcal{E}\otimes \wedge\mathcal{E}^{*}).
\end{align*}

Different choices of a connection $\nabla$ on $E \to S$ and of a degree $+1$ element $\Omega$ of $(BFV(S),[-,-]_{BFV})$
satisfying
\begin{enumerate}
\item the lowest order term of $\Omega$ is given by the tautological Section $\Omega_0$ of $\mathcal{E}\to E$ and
\item $[\Omega,\Omega]^{\nabla}_{BFV}=0$,
\end{enumerate}
lead to isomorphic
differential graded Poisson algebras 
\begin{align*}
(BFV(E),[\Omega,\cdot]^{\nabla}_{BFV},[\cdot,\cdot]^{\nabla}_{BFV}).
\end{align*}
\end{Theorem}

\begin{proof}
Pick two connections $\nabla_{0}$ and $\nabla_{1}$ on $E \to S$ and consider the two associated
graded Poisson algebras $(BFV(E),[\cdot,\cdot]^{0}_{BFV})$ and $(BFV(E),[\cdot,\cdot]^{1}_{BFV})$, respectively. 
By Corollary \ref{cor:BFV-bracket_homotopy} there is an isomorphism of graded Poisson algebras
\begin{align*}
\gamma: (BFV(E),[\cdot,\cdot]^{0}_{BFV})\xrightarrow{\cong} (BFV(E),[\cdot,\cdot]^{1}_{BFV}).
\end{align*}
Moreover
the induced automorphism of $\mathcal{C}^{\infty}(E)$ is the identity.

Assume that $\Omega$ and $\tilde{\Omega}$ are two BFV-charges of $(BFV(E),[\cdot,\cdot]^{0}_{BFV})$ and $(BFV(E),[\cdot,\cdot]^{1}_{BFV}$),
respectively. Applying the automorphism $\gamma$
to $\Omega$ yields another element of $(BFV(E),[\cdot,\cdot]^{1}_{BFV})$, which can be checked to be a BFV-charge again.
By Proposition \ref{prop:charge} this implies that there is
an inner automorphism $\beta$ of $(BFV(E),[\cdot,\cdot]^{1}_{BFV})$ which maps $\gamma(\Omega)$ to $\tilde{\Omega}$.

Hence
\begin{align*}
\beta\circ \gamma: (BFV(E),[\cdot,\cdot]^{0}_{BFV})\xrightarrow{\cong}(BFV(E),[\cdot,\cdot]^{1}_{BFV})
\end{align*}
is an isomorphism of graded Poisson algebras which maps $\Omega$ to $\tilde{\Omega}$.
\end{proof}

\section{Choice of Tubular Neighbourhood}\label{s:tubularneighbourhood}
Let $S$ be a coisotropic submanifold of a smooth, finite dimensional Poisson manifold $(M,\Pi)$.
Throughout this Section, $E$ denotes the normal bundle of $S$ inside $M$. 
As explained in subsection \ref{ss:BFV}, the first step in the construction of the BFV-complex for $S$ inside $(M,\Pi)$
is the choice of an embedding $\psi: E \hookrightarrow M$. 
Such an embedding equips $E$ with a Poisson
bivector field $\Pi_\psi$, which is used to construct the BFV-bracket on the ghost/ghost-momentum bundle,
see Subsection \ref{ss:BFV}.

Let us first consider the case where the embedding is changed by composition with a linear automorphism
of the normal bundle $E$:

\begin{Lemma}\label{lemma:linear_automorphisms}
Let 
\begin{align*}
(BFV(E),[\Omega,\cdot]_{BFV},[\cdot,\cdot]_{BFV})
\end{align*}
be a BFV-complex corresponding to some choice
of tubular neighbourhood $\psi: E \hookrightarrow M$, while
\begin{align*}
(BFV(E),[\Omega^{g},\cdot]^{g}_{BFV},[\cdot,\cdot]_{BFV}^{g})
\end{align*}
is a BFV-complex corresponding to the embedding
$\psi \circ g: E \hookrightarrow M$, where $g: E \to E$ is a vector bundle isomorphism covering the identity.

Then there is an isomorphism of graded Poisson algebras
\begin{align*}
(BFV(E),[\cdot,\cdot]_{BFV}) \to (BFV(E),[\cdot,\cdot]^g_{BFV})
\end{align*}
which maps $\Omega$ to $\Omega^{g}$.
\end{Lemma}

\begin{proof}
Let $\Pi$ / $\Pi^{g}$ be the Poisson bivector field on $E$ obtained from $\psi: E \hookrightarrow M$ / $\psi \circ g: E \hookrightarrow M$,
respectively. Clearly $\Pi^{g}=(g)_{*}(\Pi)$.

Choose some connection $\nabla$ of $E$, which is used to construct the $L_{\infty}$ quasi-isomorphism
\begin{align*}
\mathcal{L}: (\mathcal{V}(E)[1],[\cdot,\cdot]_{SN}) \leadsto (\mathcal{V}(\mathcal{E}^{*}[1]\oplus \mathcal{E}[-1])[1],[-,-]_{SN},[G,-]_{SN}).
\end{align*}
Plugging in $\Pi$ results into the BFV-bracket $[\cdot,\cdot]_{BFV}$.
On the other hand, we can use $\nabla^{g}:=(g^{-1})^{*}\nabla$ to construct another $L_{\infty}$ quasi-isomorphism $\mathcal{L}^{g}$.
Plugging in $\Pi^{g}$ results into another BFV-bracket $[\cdot,\cdot]^{g}_{BFV}$.

We claim that $[\cdot,\cdot]_{BFV}$ and $[\cdot,\cdot]^{g}_{BFV}$ are isomorphic graded Poisson brackets.
First, observe that the isomorphism $g: E \to E$ lifts to an vector bundle isomorphism
\begin{align*}
\xymatrix{
\mathcal{E} \ar[r]^{\hat{g}} \ar[d]& \mathcal{E} \ar[d]\\
E \ar[r]^{g} & E,
}
\end{align*} 
such that the tautological section gets mapped to itself under $(\hat{g})^{*}$.
We denote the induced automorphism of $E^{*}[1]\oplus E[-1]$ by $\hat{g}$ as well.

By naturality of the pull back of connections, we obtain the commutative diagram
\begin{align*}
\xymatrix{
\mathcal{V}(E) \ar[rr]^(0.4){\iota_{\nabla}} \ar[d]^{(g)_{*}}&& \mathcal{V}(E^{*}[1]\oplus E[-1]) \ar[d]^{(\hat{g})_*} \\
\mathcal{V}(E) \ar[rr]^(0.4){\iota_{\nabla^{g}}} && \mathcal{V}(E^{*}[1]\oplus E[-1]),
}
\end{align*}
where $\iota_{\nabla}$ ($\iota_{\nabla^{g}}$) is the horizontal lift induced by $\nabla$ ($\nabla^{g}$).
Using this together with the explicit description of the $L_{\infty}$ quasi-isomorphism
$\mathcal{L}$ from Proposition \ref{prop:lift} contained in \cite{Schaetz}, or in the proof of Proposition \ref{prop:lift_homotopy}, one concludes that
\begin{align*}
(\mathcal{L}^{g})_k = (\hat{g})_{*} \circ (\mathcal{L})_k \circ \left( (g)_{*}^{-1}\otimes \cdots \otimes (g)_{*}^{-1}\right). 
\end{align*}
Here, $(\mathcal{L})_k$ denotes the $k$th structure map of the $L_{\infty}$ quasi-isomorphism $\mathcal{L}$.

This immediately implies that $\hat{g}$ induces an isomorphism between $[\cdot,\cdot]_{BFV}$ and $[\cdot,\cdot]^{g}_{BFV}$,
respectively. Moreover, since $\hat{g}$ maps the tautological section to itself, it maps any BFV-charge to another one.

Finally, Theorem \ref{theoremA} implies the statement of Lemma \ref{lemma:linear_automorphisms}.
\end{proof}

In general, a different choice of embedding can cause drastic changes in the associated BFV-complexes.
Consider $S=\{0\}$ inside $M=\mathbb{R}^{2}$ equipped with the smooth Poisson bivector field 
\begin{align*}
\Pi(x,y):=\begin{cases} 0 & x^2+y^2\le 4 \\ 
\exp{(-\frac{1}{x^{2}+y^{2}-4})\frac{\partial}{\partial x}\wedge\frac{\partial}{\partial y}} & x^2+y^2\ge 4\end{cases}.
\end{align*}
Let $\psi_{0}$ be the embedding of
$E\cong \mathbb{R}^{2}$ into $\mathbb{R}^{2}$ given by the identity and $\psi_{1}$ the embedding given by
\begin{align*}
(x,y) \mapsto \frac{1}{\sqrt{1+x^{2}+y^{2}}}(x,y).
\end{align*}
The image of $\psi_{1}$ is contained in the disk of radius $1$. Hence $\Pi_{\psi_{1}}$ vanishes
identically whereas $\Pi_{\psi_{0}}$ does not.

The ghost/ghost-momentum bundle $\mathcal{E}^{*}[1]\oplus \mathcal{E}[-1]$
is of the very simple form
\begin{align*} 
\mathbb{R}^{2}\times \left( (\mathbb{R}^{2})^{*}[1]\oplus \mathbb{R}^{2}[-1]\right) \to \mathbb{R}^{2}. 
\end{align*}
Denote the Poisson bivector field coming from the natural pairing between $(\mathbb{R}^{2})^{*}[1]$ and $\mathbb{R}^{2}[-1]$ by $G$. 
We choose the standard  flat connection on the bundle $\mathbb{R}^{2}\to 0$.
Then the Poisson bivector fields for the BFV-brackets $[\cdot,\cdot]^{0}_{BFV}$ and $[\cdot,\cdot]^{1}_{BFV}$ are simply given by the sums 
$G+\Pi_{\psi_{0}}$ and $G+\Pi_{\psi_{1}}$,
respectively. 

Any isomorphism of graded Poisson algebras between
$(BFV(E),[\cdot,\cdot]^{0}_{BFV})$ and $(BFV(E),[\cdot,\cdot]_{BFV}^{1})$ yields an induced isomorphism of Poisson algebras between
$(\mathcal{C}^{\infty}(\mathbb{R}^{2}),\{\cdot,\cdot\}_{\Pi_{\psi_0}})$ and 
$(\mathcal{C}^{\infty}(\mathbb{R}^{2}),\{\cdot,\cdot\}_{\Pi_{\psi_0}})$.
Since $\Pi_{\psi_1}$ vanishes, the induced automorphism would have to map something non-vanishing to $0$, which
is a contradiction. Hence there is no isomorphism of graded Poisson algebras between $(BFV(E),[\cdot,\cdot]^{0}_{BFV})$
and $(BFV(E),[\cdot,\cdot]^{1}_{BFV})$.

Although different choices of embeddings can lead to differential graded Poisson algebras that are not isomorphic,
it is always possible to find appropriate ``restrictions'' of the BFV-complexes such that the corresponding differential
graded Poisson algebras are isomorphic. To this end we define

\begin{Definition}\label{d:restrictedBFV}
Let $E$ be a finite rank vector bundle over a smooth manifold $S$. Assume $E$ is equipped with a Poisson bivector field $\Pi$
such that $S$ is a coisotropic submanifold of $E$.
Moreover let $(BFV(E),D_{BFV},[\cdot,\cdot]_{BFV})$ be a BFV-complex for $S$ in $(E,\Pi)$ 
and $U$ an open neighbourhood of $S$ inside $E$.

Then the {\em restriction} of the BFV-complex on $U$ is the
differential graded Poisson algebra 
\begin{align*}
(BFV^{U}(E),D^{U}_{BFV}(\cdot)=[\Omega^{U},\cdot]^{U}_{BFV},[\cdot,\cdot]^{U}_{BFV})
\end{align*}
given by the following data:
\begin{itemize}
\item[(a)] $BFV^{U}(E)$ is the space of smooth functions on the graded vector bundle
$(\mathcal{E}^{*}[1]\oplus\mathcal{E}[-1])|_{U}$ fitting into the following Cartesian square:
$$
\xymatrix{
\mathcal{E}^{*}[1]\oplus\mathcal{E}[-1]|_{U} \ar[r] \ar[d] & \mathcal{E}^{*}[1] \oplus \mathcal{E}[-1] \ar[d] &\\
U \ar[r] & E.}
$$
\item[(b)] $BFV^{U}(E)$ inherits a graded Poisson bracket $[\cdot,\cdot]^{U}_{BFV}$ from $BFV(E)$: 
one restricts the Poisson bivector field
corresponding to $[\cdot,\cdot]_{BFV}$ to the graded submanifold $(\mathcal{E}^{*}[1]\oplus\mathcal{E}[-1])|_{U}$ of
$\mathcal{E}^{*}[1]\oplus\mathcal{E}[-1]$.
\item[(c)] An element $\Omega^{U}$ of $BFV^{U}(E)$ is called a {\em restricted BFV-charge} if
it is of degree $+1$, $[\Omega^{U},\Omega^{U}]^{U}_{BFV}=0$ holds and the component of $\Omega^{U}$ in $\Gamma(\mathcal{E}|_{U})$ is equal to
the restriction of the tautological section $\Omega_0 \in \Gamma(\mathcal{E})$ to $U$.
\end{itemize}
\end{Definition}

\begin{Proposition}\label{prop:central}
Let $S$ be a coisotropic submanifold of a smooth, finite dimensional Poisson manifold $(M,\Pi)$.
Denote the normal bundle of $S$ by $E$ and fix a connection $\nabla$ on $E$.
Moreover let $\psi_{0}$ and $\psi_{1}$ be two embeddings of $E$ into $M$ as tubular neighbourhoods of $S$.

Using these data one constructs two graded Poisson algebra structures on $BFV(E)$ following subsection \ref{ss:BFV} 
(in particular one applies Proposition \ref{prop:lift}). Denote the two corresponding graded Poisson brackets
by $[\cdot,\cdot]^{0}_{BFV}$ and $[\cdot,\cdot]^{1}_{BFV}$ respectively.

Then there are two open neighbourhoods $A_0$ and $A_{1}$ of $S$ in $E$ such that an isomorphism of graded Poisson algebras
\begin{align*}
(BFV^{A_{0}}(E),[\cdot,\cdot]^{0,A_{0}}_{BFV}) \xrightarrow{\cong} (BFV^{A_{1}}(E),[\cdot,\cdot]^{1,A_{1}}_{BFV})
\end{align*}
exists.
\end{Proposition}

\begin{proof}
We make use of the fact that any two embeddings of $E$ as a tubular neighbourhood are homotopic up to inner automorphisms of $E$,
i.e. given two embeddings $\psi$ and $\phi$ of $E$ into $M$ as a tubular neighbourhood, one can find
\begin{itemize}
\item a vector bundle isomorphism $g$ of $E$ and
\item a smooth map $F: E\times I \to M$
\end{itemize}
satisfying
\begin{itemize}
\item $F|_{E\times \{0\}}=\psi$ and $F|_{E\times \{1\}}=\phi \circ g$,
\item $\psi_{s}:=F|_{E\times \{s\}}: E\to M$ is an embedding for all $s\in I$ and
\item $\psi_{s}|_{S}=id_{S}$ for all $s \in I$.
\end{itemize}
The construction of $F$ can be found in \cite{Hirsch} for instance.

Since vector bundle automorphisms of $E$ yield isomorphic BVF-complexes by Lemma \ref{lemma:linear_automorphisms}, we can assume without loss 
of generality that the two embeddings $\psi:=\psi_0$ and $\phi=:\psi_1$ are homotopic (i.e. $g=id$).

Denote the images of $\psi_{s}$ by $V_s$. Since $\psi_s$ is an embedding of a manifold of the same dimension as $M$, the image
$V_{s}$ is an open subset of $M$. Moreover $S \subset V_{s}$ holds for arbitrary $s\in I$, i.e. $V_{s}$ is an open neighbourhood
of $S$ in $M$. Because $F$ is continuous, one can find an open neigbourhood $V$ of $S$ in $M$ which is contained in $\bigcap_{s\in I}V_{s}$.

One defines $\hat{F}: E\times I \to M\times I$, $(e,t)\mapsto (F(e,t),t)$ and checks that
$\hat{F}$ is an embedding, hence its image is a submanifold $W$ of $M\times I$ and $\hat{F}$ is a diffeomorphism between
$E\times I$ and $W$. Consider the restriction of $\hat{F}^{-1}: W \xrightarrow{\cong} E\times I$ to $V \times I$ which we denote by $G$.
If one restricts $G$ to ``slices'' of the form $V\times \{s\}$ one obtains $\psi_{s}^{-1}|_{V}$. The images
of $\psi_{s}^{-1}|_{V}$ are denoted by $W_{s}$. By continuity of $G$
there is an open neighbourhood $W$ of $S$ in $E$ which is contained in
$\bigcap_{s\in I}W_{s}$.

We define the following one-parameter family of local diffeomorphisms of $E$:
\begin{align*}
\phi_{s}: W_{0} \xrightarrow{\psi_{0}|_{W_{0}}} V \xrightarrow{(\psi_{s}|_{V})^{-1}} W_{s}.
\end{align*}
Moreover $E$ inherits a one-parameter family of Poisson bivector fields defined by
 $\Pi_{s}:=(\psi_{s}|_{V_s}^{-1})_{*}(\Pi|_{V_{s}})$.
The restriction $\Pi_{s}|_{W_{s}}$ is equal to $(\psi_{s}|_{V}^{-1})_{*}(\Pi|_{V})$. Consequently
\begin{align}\label{equ:flow}
\Pi_{s}|_{W_{s}}=(\phi_{s})_{*}(\Pi_{0}|_{W_0})
\end{align}
holds for all $s\in I$.

Differentiating $\phi_{s}$ yields a smooth one-parameter family of local vector fields $(Y_{s})_{s\in I}$ on $E$.
By \eqref{equ:flow} the smooth one-parameter family
\begin{align*}
\Pi_{t}|_{W} - Y_{t}|_{W}dt
\end{align*}
is a MC-element of $(\mathcal{V}(W)[1]\otimes \Omega(I),d_{DR},[\cdot,\cdot]_{SN})$.

The $L_{\infty}$ quasi-isomorphism 
\begin{align*}
\mathcal{L}_{\nabla}: (\mathcal{V}(E)[1],[\cdot,\cdot]_{SN})\leadsto (\mathcal{V}(\mathcal{E}^{*}[1]\oplus \mathcal{E}[-1])[1],[G,\cdot]_{SN},[\cdot,\cdot]_{SN})
\end{align*}
from Proposition
\ref{prop:lift} restricts to
an $L_{\infty}$ quasi-isomorphism 
\begin{align*}
\mathcal{L}_{\nabla}|_{W}: (\mathcal{V}(W)[1],[\cdot,\cdot]_{SN})\leadsto (\mathcal{V}((\mathcal{E}^{*}[1]\oplus \mathcal{E}[-1])|_{W})[1],[G,\cdot]_{SN},[\cdot,\cdot]_{SN}).
\end{align*}
Hence we obtain an $L_{\infty}$ quasi-isomorphism
\begin{eqnarray*}
\mathcal{L}_{\nabla}|_{W}\otimes id:&& (\mathcal{V}(W)[1]\otimes \Omega(I),d_{DR},[\cdot,\cdot]_{SN})\leadsto \\
&&(\mathcal{V}(\mathcal{E}^{*}[1]\oplus \mathcal{E}[-1]|_{W})[1]\otimes \Omega(I),d_{DR}+[G,\cdot]_{SN},[\cdot,\cdot]_{SN}).
\end{eqnarray*}
Applying $\mathcal{L}_{\nabla}|_{W}\otimes id$ to the MC-element $\Pi_{t}|_{W} - Y_{t}|_{W}dt$ and adding $G$
yields a MC-element $\hat{\Pi}_t - \hat{Y}_{t}dt$ of
$(\mathcal{V}(\mathcal{E}^{*}[1]\oplus \mathcal{E}[-1]|_{W})[1]\otimes \Omega(I),d_{DR},[\cdot,\cdot]_{SN})$.

It is straightforward to check that $\hat{\Pi}_s$ is the restriction of $\mathcal{L}_{\nabla}(\sum_{k\ge1}\frac{1}{k!}\Pi_{s}^{\otimes k})$
to $W$ and that $\hat{Y}_s$ is the sum of the horizontal lift $\iota_{\nabla}(Y_s)$ of $Y_{s}$ with respect to $\nabla$ restricted to $W$ plus a 
part in $\mathcal{V}^{(1,1)}(\mathcal{E}^{*}[1]\oplus \mathcal{E}[-1])$ (that acts as a nilpotent derivation).

Using parallel transport with respect to $\nabla$,
$(\iota_{\nabla}(Y_t))_{t\in I}$ can be integrated to a one-parameter family of
vector bundle automorphisms 
\begin{align*}
\hat{\phi}_{s}: \mathcal{E}^{*}[1]\oplus \mathcal{E}[-1]|_{W_{0}} \to \mathcal{E}^{*}[1]\oplus \mathcal{E}[-1]|_{W_{s}}
\end{align*}
covering $\phi_s: W_{0} \to W_{s}$ for arbitrary $s\in I$.
Similar to the construction of $V$ and $W$ one finds an open neighbourhood $A_0$ of $S$ in $W$ such that
$\phi_{t}|_{A_0}: A_0 \xrightarrow{\cong} A_{t}$ with $\bigcup_{s\in I} A_s \subset W$. So
the restriction of $\hat{\phi}_s$ to $\mathcal{E}^{*}[1]\oplus \mathcal{E}[-1]|_{A_0}$ has image
$\mathcal{E}^{*}[1]\oplus \mathcal{E}[-1]|_{A_s}$ which is a submanifold of
$\mathcal{E}^{*}[1]\oplus \mathcal{E}[-1]|_{W}$ for arbitrary $s\in I$.

Hence the one-parameter family of local vector fields 
\begin{align*}
(\iota_{\nabla}(Y_{t})|_{(\mathcal{E}^{*}[1]\oplus \mathcal{E}[-1])|_{A_t}})_{t\in I}
\end{align*}
can be uniquely integrated to a one-parameter family of local diffeomorphisms $(\hat{\phi}_{t})_{t\in I}$ and
consequently the one-parameter family of local vector fields
$(\hat{Y}_t|_{A_t})_{t\in I}$
can be uniquely integrated to a one-parameter family of local diffeomorphisms which we denote by
\begin{align*}
\varphi_{s}: (\mathcal{E}^{*}[1]\oplus \mathcal{E}[-1])|_{A_{0}} \to (\mathcal{E}^{*}[1]\oplus \mathcal{E}[-1])|_{A_{s}}
\end{align*}
for $s\in I$.

Applying Lemma \ref{l:isomorphism} shows that $\hat{\Pi}_{s}|_{A_s}=(\varphi_{s})_{*}(\hat{\Pi}_{0}|_{A_0})$ holds for all $s\in I$. Hence
\begin{align*}
(\varphi_{1})_{*}: \mathcal{C}^{\infty}(\mathcal{E}^{*}[1]\oplus \mathcal{E}[-1]|_{A_{0}}) \to
\mathcal{C}^{\infty}(\mathcal{E}^{*}[1]\oplus \mathcal{E}[-1]|_{A_{1}})
\end{align*}
is an isomorphism of Poisson algebras. 
\end{proof}

\begin{Theorem}\label{thm:central}
Let $S$ be a coisotropic submanifold of a smooth, finite dimensional Poisson manifold $(M,\Pi)$.
Suppose $(BFV(E),D^{0}_{BFV},[\cdot,\cdot]^{0}_{BFV})$ and
$(BFV(E),D^{1}_{BFV},[\cdot,\cdot]^{1}_{BFV})$ are two BFV-complexes
constructed with help of two arbitrary embeddings of $E$ into $M$, two arbitrary connections on $E\to S$
and two arbitrary BFV-charges.

Then there are two open neighbourhoods $B_0$ and $B_1$ of $S$ in $E$ such that an isomorphism of
differential graded Poisson algebras
\begin{align*}
(BFV^{B_{0}}(E),D_{BFV}^{0,B_{0}},[\cdot,\cdot]^{0,B_{0}}_{BFV}) \xrightarrow{\cong} (BFV^{B_{1}}(E),D_{BFV}^{1,B_{1}}[\cdot,\cdot]^{1,B_{1}}_{BFV})
\end{align*}
exists.
\end{Theorem}

\begin{proof}
By Theorem \ref{theoremA} we can assume without loss of generality that the two chosen connections coincide.
Furthermore it suffices to prove that there is an isomorphism of graded Poisson algebras
from some restriction of $(BFV(E),[\cdot,\cdot]^{0}_{BFV})$ to some restriction of
$(BFV(E),[\cdot,\cdot]^{0}_{BFV})$ which maps a restricted BFV-charge to another restricted BFV-charge.
This is a consequence of the fact that Theorem \ref{theoremA} holds also in the restricted setting as long as the 
open neighbourhood $U$ of $S$ in $E$, to which we restrict, is contractible to $S$ along the fibres of $E$.

By Lemma \ref{lemma:linear_automorphisms}, we may assume without loss of generality that the two
embeddings under consideration are homotopic. Hence there is
a smooth one-parameter family of isomorphisms of graded Poisson algebras
\begin{align*}
(\varphi_{s})_{*}: (BFV^{A_0}(E),[\cdot,\cdot]^{0,A_0}_{BFV}) \to (BFV^{A_s}(E),[\cdot,\cdot]^{s,A_s}_{BFV}),
\end{align*}
which we constructed in the proof of Proposition \ref{prop:central}.
The smoothness of this family and the fact that the zero section $S$ is fixed under $(\varphi_{s})_{s\in I}$ 
imply that there is a open neighbourhood $A$ of $S$ in $E$ satisfying
$A\subset \bigcap_{s\in I} A_s$.

Fix a restricted BFV-charge $\Omega$ of $(BFV^{A_{0}}(E),[\cdot,\cdot]^{0,A_0}_{BFV})$. The restriction of
\begin{align*}
(\Omega(t):=(\varphi_{t})_{*}(\Omega))_{t\in I}
\end{align*}
to $A$ yields a smooth one-parameter family of sections of $\bigwedge \mathcal{E} \otimes
\bigwedge \mathcal{E}^{*}|_{A}$. Although $[\Omega(s)|_A,\Omega(s)|_A]_{BFV}^{s,A}=0$ holds for all $s\in I$,
$\Omega(s)|_{A}$ is in general not a BFV-charge since its component in $\Gamma(\mathcal{E}|_{W})$ is
$\Omega_0(s):=(\varphi_{s})_{*}(\Omega_0)$ which does not need to be equal to $\Omega_0$ as required -- see Definition \ref{d:restrictedBFV}.
In particular $\Omega(1)$ might not be a restricted BFV-charge of $(BFV(E),[\cdot,\cdot]^{1}_{BFV})$.
However we will show that $\Omega(1)$ can be ``gauged'' to a BFV-charge in the remainder of the proof.

We have to recall some of the ingredients involved in the proof of Proposition \ref{prop:charge}:
The first observation is that $\delta:=[\Omega_0,\cdot]_{G}$ is a differential. Here $\Omega_0$ denotes
the tautological section of $\mathcal{E}\to E$, $G$ is the Poisson bivector field associated to the
fibre pairing between $\mathcal{E}$ and $\mathcal{E}^{*}$, and $[\cdot,\cdot]_{G}$ denotes the graded Poisson bracket
on $BFV(E)$ corresponding to $G$. Second it is possible to construct a homotopy $h$ for $\delta$, i.e. a degree $-1$
map satisfying
\begin{align}\label{homotopy}
\delta \circ h + h \circ \delta = id - i\circ pr
\end{align}
where $i$ is an embedding of the cohomology of $\delta$ into $BFV(E)$ and $pr$ is a projection from $BFV(E)$ onto cohomology.
We remark that $h$ does not restrict to arbitrary open neighbourhoods of $S$ in $E$. However one can check that it does restrict
to open neighbourhoods that can be contracted to $S$ along the fibres of $E$.
Without loss of generality we can assume that $A$ has this property.

We are interested in the smooth one-parameter family 
\begin{align*}
h(\Omega_{0}(s))\in \Gamma(\mathcal{E}\otimes \mathcal{E}^{*}|_A)\cong \Gamma(\End(\mathcal{E}|_{A}))
\end{align*}
with $s \in I$.
Since $\Omega_0$ intersects the zero section of $\mathcal{E} \to E$ transversally at $S$, so does
$\Omega_0(s)$ for arbitrary $s\in I$. This implies 1.) the evaluation of $\Omega_0(s)$ at $S$ is zero and
2.) $h(\Omega_0(s))|_{S} \in \Gamma(\mathcal{E}\otimes \mathcal{E}^{*}|_{S})$ is fibrewise invertible, i.e.
it is an element of $\Gamma(GL(\mathcal{E}|_{S}))$.

For any $s \in I$ we have $\delta(\Omega_0(s))=[\Omega_0,\Omega_0(s)]_G=0$ since both $\Omega_0$ and
$\Omega_0(s)$ are sections of $\mathcal{E}|_{A}$ and $G$ is the Poisson bivector given by contraction between
$\mathcal{E}$ and $\mathcal{E}^{*}$. Moreover $(i\circ pr)(\Omega_0(s))=0$ since the projection $pr$ involves
evaluation of the section at $S$, where $\Omega_0(s)$ vanishes.
Consequently \eqref{homotopy} reduces to $\delta(h(\Omega_0(s)))=\Omega_0(s)$ for all $s\in I$.
However this means that if we interpret $h(\Omega_0(s))$ as a fibrewise endomorphism of $\mathcal{E}|_{A}$
the image of $\Omega_0$ under $-h(\Omega_0(s))$ is $\Omega_0(s)$.

We define $M_s:=-h(\Omega_0(s))$ -- as already observed, $(M_t)_{t\in I}$ is a smooth one-parameter family of
sections of $\End(\mathcal{E}|_{A})$ and the restriction to $S$ is a smooth one-parameter family of
$GL(\mathcal{E}|_{S})$. By smoothness of the one-parameter family it is possible to find an open neighbourhood
$B$ of $S$ in $E$ such that the restriction of $(M_t)_{t\in I}$ to $B$ is always fibrewise invertible. Since
$M_0=id|_{A}$ we know that $(M_t|_{B})_{t\in I}$ is a smooth one-parameter family of sections in $GL_{+}(\mathcal{E}|_{B})$, i.e.
fibrewise invertible automorphisms of $E|_{B}$ with positive determinante.
In particular $M_{1} \in \Gamma(GL_{+}(\mathcal{E}|_{B}))$. 

Consider the smooth one-parameter family $(m_{t})_{t\in I}$ of sections of $\End(\mathcal{E}|_{B})$ given by
\begin{align*}
m_{t}:=-M_{t}^{-1}\circ \left(\frac{d}{dt}M_{t}\right).
\end{align*}
It integrates to a smooth one-parameter family of sections of $GL_{+}(\mathcal{E}|_{B})$ that coincides with $(M_{t})_{t\in [0,1]}$.
The adjoint action of $m_t$ on $(BFV^{B}(E),[\cdot,\cdot]_{BFV}^{1,B})$ can be integrated to an automorphism of $(BFV^{B}(E),[\cdot,\cdot]_{BFV}^{1,B})$ and this automorphism maps the restriction of $\Omega_0(1)$ to $B$ to 
the restriction of $\Omega_0$ to $B$. Hence $(\exp(m)\circ (\varphi_{1})_{*})$ maps the restricted BFV-charge $\Omega$
to another restricted BFV-charge of $(BFV^{B}(E),[\cdot,\cdot]_{BFV}^{1,B})$.
\end{proof}

\begin{Definition}\label{d:germBFV}
Let $(BFV(E),D_{BFV},[\cdot,\cdot]_{BFV})$ be a BFV-complex associated to a coisotropic submanifold $S$
of a smooth Poisson manifold $(M,\Pi)$. We define a differential graded Poisson algebra
$(BFV^{\mathfrak{g}}(E),D^{\mathfrak{g}}_{BFV},[\cdot,\cdot]^{\mathfrak{g}}_{BFV})$ as follows:
\begin{itemize}
\item[(a)] $BFV^{\mathfrak{g}}(E)$ is the algebra of equivalence classes of elements of $BFV(E)$ under the equivalence
relation: $f\sim g :\Leftrightarrow$ there is a open neighbourhood $U$ of $S$ in $E$ such that $f|_{U}=g|_{U}$.
\item[(b)] $D^{\mathfrak{g}}_{BFV}([\cdot]):=[D_{BFV}(\cdot)]$ where $[\cdot]$ denotes the equivalence class of $\cdot$ under $\sim$.
\item[(c)] $[[\cdot],[\cdot]]^{\mathfrak{g}}_{BFV}:=[[\cdot,\cdot]_{BFV}]$.
\end{itemize}
\end{Definition}

Given a differential graded Poisson algebra with unit $(A,\wedge,d,[\cdot,\cdot])$ we define the corresponding 
{\em abstract differential graded Poisson algebra with unit} $[(A,\wedge,d,[\cdot,\cdot])]$ to be the isomorphism class of
$(A,\wedge,d,[\cdot,\cdot])$ in the category of differential graded Poisson algebras with unit. In particular
$[(A,\wedge,d,[\cdot,\cdot])]$ is a object in the category of differential graded Poisson algebras
with unit up to isomorphisms.

Theorem \ref{thm:central} immediately implies

\begin{Corollary}\label{cor:central}
Consider a coisotropic submanifold $S$ of a smooth, finite dimensional Poisson manifold $(M,\Pi)$ and 
let $(BFV(E),D_{BFV},[\cdot,\cdot]_{BFV})$ be a BFV-complex associated to $S$ inside $(M,\Pi)$.

The abstract differential graded Poisson algebra
\begin{align*}
[(BFV^{\mathfrak{g}}(E),D^{\mathfrak{g}}_{BFV},[\cdot,\cdot]^{\mathfrak{g}}_{BFV})]
\end{align*}
is independent of the specific choice of a BFV-complex and hence is an invariant of
$S$ as a coisotropic submanifold of $(M,\Pi)$.
\end{Corollary}

\end{document}